\documentclass[notitlepage,10pt]{article}
\usepackage{amssymb,mathrsfs,amsmath,enumitem}
\catcode`\@=11
\@addtoreset{equation}{section}

\catcode`\@=12

\newtheorem{Theorem}{Theorem}[section]
\newtheorem{Lemma}[Theorem]{Lemma}
\newtheorem{Proposition}[Theorem]{Proposition}

\newtheorem{Remark}[Theorem]{Remark}
\newtheorem{Definition}[Theorem]{Definition}

\def\R{{\mathbb R}}

\def\N{{\mathbb N}}
\def\S{{\mathbb S}}

\def\A{{\mathcal A}}
\def\B{{\mathcal B}}
\def\D{{\mathcal D}}
\def\E{{\mathcal E}}
\def\F{{\mathcal F}}

\def\Q{{\mathcal Q}}
\def\V{{\mathcal V}}

\def\Proof{\noindent\textit{Proof. }}
\def\qed{$~\square$\goodbreak \medskip}

\title{Existence of isovolumetric extremals for capillarity functionals}
\author{Paolo Caldiroli\footnote{Dipartimento di Matematica, Universit\`a di Torino, via Carlo Alberto, 10 -- 10123 Torino, Italy. Email: \tt{paolo.caldiroli@unito.it}}
, Alessandro Iacopetti\footnote{Dipartimento di Matematica, Universit\`a di Torino, via Carlo Alberto, 10 -- 10123 Torino, Italy. Email: \tt{alessandro.iacopetti@unito.it}}}

\begin{document}

\date{}
\maketitle

\begin{abstract}
\noindent
Capillarity functionals are parameter invariant functionals defined on classes of two-dimensional parametric surfaces in $\R^{3}$ as the sum of the area integral and a non homogeneous term of suitable form. Here we consider the case of a class of non homogenous terms vanishing at infinity for which the corresponding capillarity functional has no volume-constrained $\S^{2}$-type minimal surface. Using variational techniques, we prove existence of extremals characterized as saddle-type critical points. 
\smallskip

\smallskip

\end{abstract}

\section{Introduction}
Surfaces of constant mean curvature are critical points of the area functional for volume-preserving variations. They constitute a nice model for describing closed capillarity surfaces, i.e. soap bubbles, when the surface energy of the liquid is regarded as isotropic, the liquid is homogeneous and no external force is considered. In this case the surface energy is proportional to the surface area, and soap bubbles correspond to extremal solutions of the isoperimetric problem.

If external forces are taken into account, then the surface energy has to be modified in a suitable way, by considering a generalized area functional
\begin{equation}
\label{eq:f-area}
A_w(\Sigma)= \int_\Sigma w(p) \ d\Sigma\!~,
\end{equation}
where $w\colon\R^{3}\to\R$ is a regular and positive weight.

Functionals of the form (\ref{eq:f-area}) have been extensively studied from the viewpoint of geometric measure theory (as in \cite{CafDel01}, for instance). Correspondingly, in the same direction, also isoperimetric problems with weights have been recently studied, in some cases (see \cite{Morg}, \cite{MorgPrat}).

Here we are interested in investigating some issues about a class of generalized area functionals, from a different perspective, in the frame of differential geometry. With this approach we are allowed to prescribe the topological type of the surfaces we deal with. In particular, we focus on parametric surfaces of the type of the sphere. This means that we identify surfaces with (the range of) maps from $\S^{2}$ to $\R^{3}$. Moreover we consider functionals of the kind
$$
F(u)= \int_{\S^2} (1+ Q(u)\cdot \nu) \ d\Sigma~\!,
$$
where $\nu$ is the Gauss map, $d\Sigma$ is the area element of $\S^2$ induced by $u$, and $Q\colon\R^3\to\R^3$ is a prescribed smooth vector field such that
\begin{equation}
\label{eq:Q-bound}
\|Q\|_\infty < 1~\!.
\end{equation} 
These functionals are known as ``capillarity functionals'' (see \cite{HilMos}) and they can be seen as a correction of the area functional by a non homogeneous term. The bound (\ref{eq:Q-bound}) is a sufficient (and necessary) condition in order that an isoperimetric inequality for capillarity functionals holds true. We are interested in looking for critical points for these kind of functionals in the Sobolev space $H^{1}(\S^{2},\R^{3})$, for volume-preserving variations, assuming that the non homogeneous term vanishes at infinity, namely
\begin{equation}
\label{eq:Q-vanish}
Q(p)\to 0\quad\text{as }|p|\to\infty~\!.
\end{equation} 
Actually, we can state the precise assumptions just on the scalar field $K=\mathrm{div}~\!Q$, because capillarity functionals depend on the vector field $Q$ only by its divergence.

In fact, the datum of our problem is a regular enough, scalar field $K\colon\R^{3}\to\R$ satisfying: 
\begin{itemize}[leftmargin=25pt]
\item[$(K_{1})$]~~$\sup_{p\in\R^{3}}|K(p)p|=:k_{0}< 2$ for every $p\in\R^{3}$.
\item[$(K_{2})$]~~$K(p)p\to 0$ as $|p|\to\infty$. 
\end{itemize}
Then it is possible to construct a vector field $Q_{K}\in C^{1}(\R^{3},\R^{3})$ such that $\mathrm{div}~\!Q_{K}=K$ on $\R^{3}$ and satisfying (\ref{eq:Q-bound}) and (\ref{eq:Q-vanish}) which are direct consequences of $(K_{1})$ and $(K_{2})$, respectively (see Remark \ref{R:K-volume}). For this reason, assumptions $(K_{1})$ and $(K_{2})$ seem to be reasonably natural to deal with situations with non homogeneous terms vanishing at infinity.

In general, even if the non homogeneous term vanishes at infinity, its presence in the capillarity functional has important consequences on the issue of the existence of extremals for the corresponding isoperimetric inequality. In \cite{Cald2015} one can find some results concerning both existence and  non-existence of critical points corresponding to minima for the isoperimetric problems
\begin{equation}
\label{eq:SKt}
\begin{split}
\mathcal{S}_{K}(t):=\inf\left\{\F_{K}(u)~|~u\in H^{1}(\S^{2},\R^{3}),~\V(u)=t\right\}\\
\text{where}\quad\F_{K}(u):= \int_{\S^2} (1+ Q_K(u) \cdot \nu) \ d\Sigma
\end{split}
\end{equation}
and $\V(u)$ is the algebraic volume functional, defined as the unique continuous extension to $H^{1}(\S^{2},\R^{3})$ of the integral functional 
$$
\V(u)=\frac{1}{3}\int_{\S^{2}} u \cdot \nu ~\!d\Sigma\quad\text{for }u\in H^{1}(\S^{2},\R^{3})\cap L^{\infty}.
$$ 
For future convenience, let us state some results proved in \cite{Cald2015}, about problems (\ref{eq:SKt}) with $t>0$.

\begin{Theorem}
\label{T:existence1}
Let $K\in C^{1}(\R^{3})$ satisfy $(K_{1})$--$(K_{2})$. Moreover assume that 
\begin{equation}
\label{eq:sign<0}
K(p)<0\quad\text{at some }p\in\R^{3}
\end{equation}
and that the constant $k_{0}$ appearing in $(K_{1})$ satisfies
\begin{equation}
\label{eq;restriction}
2^{2/3}(2+k_{0})<(2-k_{0})^{2}~\!.
\end{equation}
Then there exists $t_{+}>0$ such that for every $t\in(0,t_{+})$ the minimization problem defined by (\ref{eq:SKt}) admits a minimizer.
\end{Theorem}
The value $t_{+}$ can be characterized as follows
$$
t_{+}:=\sup\left\{t\ge 0~|~K\le 0\text{ and }K\not\equiv 0\text{ in some ball of radius }\sqrt[3]{3t/4\pi}\right\}.
$$
In particular $t_{+}=\infty$ if $K\le 0$ everywhere (but also if $K\le 0$ on the tail of an open cone). Other conditions on $K$, different from (\ref{eq;restriction}) and regarding the radial oscillation of $K$ are also displayed in \cite{Cald2015}. Moreover in \cite{Cald2015} it is proved that
\begin{Theorem}
\label{T:nonexistence}
Let $K\in C^{0}(\R^{3})$ satisfy $(K_{1})$--$(K_{2})$. If 
\begin{equation}
\label{eq:sign>0}
K(p)>0\quad\text{for every }p\in\R^{3},
\end{equation}
then there exists $\tau>0$ such that for every $t\in(0,\tau)$ the minimization problem defined by (\ref{eq:SKt}) has no minimizer. Moreover ${S}_{K}(t)=S t^{2/3}$ for $t\in(0,\tau)$, where $S=\sqrt[3]{36\pi}$ is the isoperimetric constant. 
\end{Theorem}

The present paper is a continuation and a completion of \cite{Cald2015}. Here we focus on the issue of existence of critical points in the case of nonexistence of minima. 



\begin{Theorem}\label{mainteo}
Let $K \in C^{1,\alpha}(\R^3)$ satisfy $(K_1)$--$(K_2)$. Moreover assume  (\ref{eq:sign>0}) and that the constant $k_{0}$ appearing in $(K_{1})$ satisfies 
\begin{equation}\label{condkzero}
k_{0}<2 (2^{1/3}-1).
\end{equation}
Then there exists a sequence $t_{n}\to 0^{+}$ such that the set of constrained critical points of $\F_K$ at volume $t_{n}$, denoted $\mathrm{Crit}_{\F_{K}}(t_{n})$, is non empty. 
\end{Theorem}





 
The proof of this result is mainly based on a min-max argument and on  degree theory, in the spirit of a procedure introduced in \cite{BahriLi} for certain semilinear elliptic equations in $\R^{N}$. 

More precisely, arguing by contradiction, if there are no volume-constrained critical points, we can construct a suitable minimax level $c$ for the functional which lies between two consecutive levels, corresponding to the the energy at infinity, i.e. the area, of one and two identical spheres at fixed volume. On the other hand, if there are no volume-constrained critical points, then constrained Palais-Smale sequences have a limit configuration made by a finite number of spheres, each one carrying the same energy. This fact comes out by some key results obtained in \cite{CaMuJFA} and \cite{CaJFA}. Hence the contradiction follows by proving the existence of volume-constrained Palais-Smale sequences at the minimax level $c$  (see Proposition \ref{propconstrminmaxPS}). 

We stress that the existence of volume-constrained Palais-Smale sequences at the minimax level $c$ is a delicate and rather technical step. In fact, in general, $\mathcal{F}_K$ is not $C^1$ and not even Gateaux differentiable. To our knowledge, a similar result is only available in the context of minimax levels for the free functionals (see \cite{MawWil89}) and only for $C^{1}$ functionals. Furthermore, a constrained version of our Proposition can be obtained through a deformation-lemma argument but it requires the functional to be of class $C^1$ and the constraint to be a Finsler manifold of class $C^{1,1}$. Instead our proof is just based on the Ekeland's variational principle (see, e.g., \cite{Eke79}) and fine estimates. 

We also point out that capillarity functionals are particularly meaningful because of their connection with the $H$-bubble problem. In fact, volume-constrained extremals parametrize $\S^{2}$-type surfaces with volume $t$ and mean curvature $H(p)=\frac{1}{2}(K(p)-\lambda)$, where $K=\mathrm{div}~\!Q$ is  prescribed, and $\lambda$ is a constant corresponding to the Lagrange multiplier due to the constraint. Differently from previous results obtained for the $H$-bubble problem, the mean curvature is prescribed up to a constant, while  in \cite{CaMuCCM}, \cite{CaMu11} the mean curvature is of the form $H(p)=\frac{1}{2}(K(p)-\lambda_0)$, where $\lambda_0$ is a given constant but no information is provided on the volume of those surfaces. In addition, it is important to note that in our paper we just assume $(K_{1})$ with \eqref{condkzero} and $(K_{2})$ (see Theorem \ref{mainteo2}), while in \cite{CaMuCCM} and \cite{CaMu11}, for analogous results one needs more restrictive assumptions, involving the radial derivative of $K$.

We also point out that even though we obtain an existence result only for a sequence $t_{n} \to 0^+$, we believe that our result is relevant in view of the techniques applied for the proof. We suspect that other methods, as the finite-dimensional reduction method, could be used to get an existence result for all $t$ in a small interval $(0,\epsilon)$. By the way, this strategy, already employed for the $H$-bubble problem (see, e.g. \cite{CaMM}, \cite{CaMuDuke}, \cite{Fel05}, \cite{MuJAM}) has not been investigated so far for the generalized isoperimetric problem.

A great part of the tools we use in the present paper is contained in \cite{Cald2015} and for the sake of convenience we recall them in Section 2. Sections 3, 4 and 5 are devoted, respectively, to the construction of the minimax scheme, to the existence of constrained Palais-Smale sequences and to the proof of Theorem \ref{mainteo}.

\section{Preliminaries}
\label{S:preliminaries}
Let us introduce the space
$$
\hat{H}^{1}(\R^2,\R^3):=\{u\in H^{1}_{loc}(\R^{2},\R^{3}); \ \int_{\R^{2}}(|\nabla u|^{2}+\mu^{2}|u|^{2})<\infty\},
$$ 
where 
\begin{equation}
\label{eq:mu2}
\mu(z)=\frac{2}{1+|z|^{2}}\quad\text{for } z=(x,y)\in\R^{2}. 
\end{equation}
For simplicity we will use the notation $\hat{H}^1$ instead of $\hat{H}^{1}(\R^2,\R^3)$.
The space $\hat{H}^{1}$ is a Hilbert space with inner product
$$
\langle u,v\rangle=\int_{\R^{2}}(\nabla u\cdot\nabla v)+\left(\frac{1}{4\pi}\int_{\R^{2}}u\mu^{2}\right)\cdot\left(\frac{1}{4\pi}\int_{\R^{2}}v\mu^{2}\right)
$$
and is isomorphic to the space $H^{1}(\S^{2},\R^{3})$. The isomorphism is given by the correspondence $\hat{H}^{1}\ni u\mapsto u\circ\phi\in H^{1}(\S^{2},\R^{3})$, where $\phi$ is the stereographic projection of $\S^{2}$ onto the compactified plane $\R^{2}\cup\{\infty\}$. As usual, we denote $\|u\|=\langle u,u\rangle^{1/2}$.

It is known that $C^{\infty}(\S^{2},\R^{3})$ is dense in $H^{1}(\S^{2},\R^{3})$ (see, e.g., \cite{Aub}, Ch.2). As a consequence, 
$\hat{C}^{\infty}:=\{u\circ\phi^{-1}~|~u\in C^{\infty}(\S^{2},\R^{3})\}$ is dense in $\hat{H}^{1}$. We point out that constant maps belong to $\hat{H}^1$, and we identify them with $\R^3$. Moreover we observe that $p+\hat{H}^{1}=\hat{H}^{1}$ for every $p\in\R^{3}$.\\ 

We recall now some important facts. Some of them are well known and classical. Others, more related to our problem, are discussed in \cite{Cald2015}. We refer to that paper for the proofs or for additional, useful bibliography.

\begin{Lemma}
\label{L:density}
The space $\R^{3}+C^{\infty}_{c}(\R^{2},\R^{3})$ is dense in $\hat{H}^{1}$. In particular, for every $u\in\hat{H}^{1}\cap L^{\infty}$ there exists a sequence $(u^{n})\subset\R^{3}+C^{\infty}_{c}(\R^{2},\R^{3})$ such that $u^{n}\to u$ in $\hat{H}^{1}$, in $L^{\infty}_{loc}$ and $\|u^{n}\|_{\infty}\le\|u\|_{\infty}$.
\end{Lemma}

Set
$$
\A(u):=\int_{\R^{2}}|u_{x}\wedge u_{y}|,~~ \D(u):=\frac{1}{2}\int_{\R^{2}}|\nabla u|^{2}\quad(u\in\hat{H}^{1})~~\text{and}~~\V(u):=\frac{1}{3}\int_{\R^{2}}u\cdot u_{x}\wedge u_{y}\quad(u\in\hat{H}^{1}\cap L^{\infty}).
$$

\begin{Lemma}
\label{L:volume}
The functional $\V$ admits a unique analytic extension on $\hat{H}^{1}$. In particular for every $u\in\hat{H}^{1}$
$$
\V'(u)[\varphi]=\int_{\R^{2}}\varphi\cdot u_{x}\wedge u_{y}\quad\forall\varphi\in \hat{H}^{1}\cap L^{\infty}
$$
and there exists a unique $v\in\hat{H}^{1}\cap L^{\infty}$ which is a (weak) solution of 
$$
\left\{\begin{array}{l}
-\Delta v=u_{x}\wedge u_{y}\quad\text{on }\R^{2}\\
\int_{\R^{2}}v\mu^{2}=0.
\end{array}\right.
$$
Moreover 
\begin{equation}
\label{eq:Wente-inequality}
\|\nabla v\|_{2}+\|v\|_{\infty}\le C\|\nabla u\|_{2}^{2}
\end{equation}
for a constant $C$ independent of $u$. 
In addition, for every $t\ne 0$ the set 
\begin{equation}
\label{eq:Mt-def}
M_{t}:=\{u\in\hat{H}^{1}~|~\V(u)=t\}
\end{equation}
is a smooth manifold and, for any fixed $u\in M_{t}$, a function $\varphi\in\hat{H}^{1}$ belongs to the tangent space to $M_{t}$ at $u$, denoted $T_{u}M_{t}$, if and only if $\V'(u)[\varphi]=0$. 
\end{Lemma}

\begin{Remark}
\label{R:volume}
The second part of Lemma \ref{L:volume} states that there exists $C>0$ such that $\|\V'(u)\|_{\hat{H}^{-1}}\le C\|\nabla u\|_{2}^{2}$ for every $u\in\hat{H}^{1}$, where $\hat{H}^{-1}$ denotes the dual of $\hat{H}^{1}$. 
\end{Remark}

\begin{Remark}
\label{R:sphere}
The mapping $\omega(z)=\left(\mu x,\mu y,1-\mu\right)$, with $\mu$ defined in (\ref{eq:mu2}), is a conformal parame\-tri\-za\-tion of the unit sphere. Indeed, it is the inverse of the stereographic projection from the North Pole.  Moreover $\A(\omega)=\D(\omega)=4\pi$ and $\V(\omega)=-\frac{4\pi}{3}$. If $\overline{p}\in\R^{3}$ and $r\in\R\setminus\{0\}$, then $u=\overline{p}+r\omega$ is a parametrization of a sphere centered at $\overline{p}$ and with radius $|r|$, Moreover $\A(u)=\D(u)=4\pi r^{2}$ and $\V(u)=-\frac{4\pi r^{3}}{3}$.
\end{Remark}

\begin{Lemma}[Isoperimetric inequality]
\label{L:isoperimetric}
It holds that
\begin{equation}
\label{eq:isoperimetric}
S|\V(u)|^{2/3}\le\A(u)\le\D(u)\quad\forall u\in \hat{H}^{1},
\end{equation}
where $S=\sqrt[3]{36\pi}$ is the best constant. Moreover any extremal function for (\ref{eq:isoperimetric}) is a conformal parametrization of a simple sphere. 
\end{Lemma}

\medskip

Fixing $K\in C^{1}(\R^{3})$ satisfying $(K_{1})$, set 
\begin{equation}\label{campodivk}
m_{K}(p):=\int^{1}_{0}K(sp)s^{2}~\!ds\quad\text{and}\quad Q_{K}(p):=m_{K}(p)p\quad\forall p\in\R^{3}
\end{equation}
and observe that $\mathrm{div}~\!Q_{K}=K$. 
Then set
$$
\Q(u):=\int_{\R^{2}}Q_{K}(u)\cdot u_{x}\wedge u_{y}\quad(u\in\hat{H}^{1}).
$$
\begin{Remark}\label{R:corresp}
We point out that under the correspondence $u\mapsto u\circ\phi$ it holds that
$$\F_{K}(u \circ \phi)= \int_{\S^2} (1+ Q_K(u\circ\phi)\cdot\nu) \ d\Sigma=\int_{\R^{2}}\left(|u_{x}\wedge u_{y}| + Q_K(u) \cdot u_{x}\wedge u_{y} \right)~\!dx~\!dy = \A(u) + \Q(u).$$
In view of this equality we can extend $\F_{K}$ to the class of non immersed surfaces. 
\end{Remark}

\begin{Remark}
\label{R:K-volume}
Using (\ref{campodivk}) one can easily check that $Q_{k}$ satisfies (\ref{eq:Q-bound}). More precisely, 
\begin{equation}
\label{eq:QK-bound}
\|Q_{K}\|_{\infty}\le\frac{k_{0}}{2}<1.
\end{equation}
Moreover the functional $\Q$ is well defined on $\hat{H}^{1}$ and 
\begin{equation}
\label{eq:VK-estimate}
|\Q(u)|\le\|Q_{K}\|_{\infty}\D(u)\quad\forall u\in\hat{H}^{1}.
\end{equation}
One can also check that $Q_{k}$ satisfies (\ref{eq:Q-vanish}). Indeed, for $|p|>R$ write 
$$
Q_{K}(p)=\frac{\hat{p}}{|p|^{2}}\int_{0}^{R}K(t\hat{p})t^{2}dt+\frac{\hat{p}}{|p|^{2}}\int_{R}^{|p|}K(t\hat{p})t^{2}dt
$$
with $\hat{p}=\frac{p}{|p|}$, and use $(K_{2})$ to conclude. 
\end{Remark}

The next result collects some useful properties of the functional $\Q$.

\begin{Lemma}
\label{L:K-volume}
Let $K\colon\R^{3}\to\R$ be a bounded continuous function. Then:
\begin{itemize}[leftmargin=22pt]
\item[(i)]
the functional $\Q$ is continuous in $\hat{H}^{1}$.
\item[(ii)]
For every $u\in\hat{H}^{1}$ and $\varphi\in\hat{H}^{1}\cap L^{\infty}$ one has 
$$
\Q(u+\varphi)-\Q(u)=\int_{\R^{2}}\left(\int_{0}^{1}K(u+r\varphi)\varphi\cdot(u_{x}+r\varphi_{x})\wedge(u_{y}+r\varphi_{y})~\!dr\right)~\!dx~\!dy.
$$
\item[(iii)]
The functional $\Q$ admits directional derivatives at every $u\in\hat{H}^{1}$ along any $\varphi\in\hat{H}^{1}\cap L^{\infty}$, given by
$$
\Q'(u)[\varphi]=\int_{\R^{2}}K(u)\varphi\cdot u_{x}\wedge u_{y}.
$$
\end{itemize}
If in addition $\sup_{p\in\R^{3}}|K(p)p|<\infty$ then for every $u\in\hat{H}^{1}$ the mapping $s\mapsto\Q(su)$ is differentiable and
$$
\frac{d}{ds}[\Q(su)]=s^{2}\int_{\R^{2}}K(su)u\cdot u_{x}\wedge u_{y}.
$$
\end{Lemma}

Now we state and prove a technical result which will be useful in the sequel.

\begin{Lemma} \label{techlemmacont}
For any $\varphi \in \R^3 + C_0^\infty(\R^2,\R^3)$ the map $u \mapsto \E^\prime(u)[\varphi]$ from $\hat{H}^1$ to $\R$ is continuous.
\end{Lemma}

\Proof
Thanks to Lemma \ref{L:K-volume} (iii) we have that for any $u \in \hat{H}^1$ and $\varphi \in \R^3 + C_0^\infty(\R^2,\R^3)$ the functional  $\E$ admits the directional derivative at $u$ along $\varphi$ and
$$
\E'(u)[\varphi]=\int_{\R^2} \nabla u \cdot \nabla \varphi + \int_{\R^{2}}K(u)\varphi\cdot u_{x}\wedge u_{y}.
$$
Since $u\mapsto \D^\prime(u)[\varphi]$ is continuous, it suffices to show that this holds also for $u\mapsto \Q'(u)[\varphi]$. Let $(u^n) \subset \hat{H}^1$ be such that $u^n \to u$ in $\hat{H}^1$. Then
\begin{equation*}
\begin{split}
| \Q'(u^n)&[\varphi]- \Q'(u)[\varphi]|\\
&=\left| \int_{\R^{2}}K(u^n)\varphi\cdot u^n_{x}\wedge u^n_{y} -  \int_{\R^{2}}K(u)\varphi\cdot u_{x}\wedge u_{y} \right|\\
&=\left| \int_{\R^{2}}K(u^n)\varphi\cdot u^n_{x}\wedge u^n_{y} -  \int_{\R^{2}} K(u)\varphi\cdot u^n_{x}\wedge u^n_{y} +  \int_{\R^{2}} K(u)\varphi\cdot u^n_{x}\wedge u^n_{y} -  \int_{\R^{2}}K(u)\varphi\cdot u_{x}\wedge u_{y} \right|\\
&\leq \int_{\R^{2}}|K(u^n)-K(u)| |\varphi| |u^n_{x}\wedge u^n_{y}| + \int_{\R^{2}}|K(u)| |\varphi| |u^n_{x}\wedge u^n_{y} - u_{x}\wedge u_{y} |=I^n_1+I^n_2.
\end{split}
\end{equation*}
Since $u^n \to u$ in $\hat{H}^1$, we get that $u^n_{x}\wedge u^n_{y} \to u_{x}\wedge u_{y}$ in $L^1(\R^2)$ and $u^n \to u$ a.e. in $\R^2$. Moreover, since $K$ is continuous and satisfies $(K_1)$, then $K$ is bounded by some positive constant $C_K$. Now assume by contradiction that $I^n_1 \not\to 0$ as $n \to\infty$. This means that there exists $\epsilon>0$ such that $|I^{n_k}_1| > \epsilon$ for some subsequence $n_k \to\infty$. But  since $u^n_{x}\wedge u^n_{y} \to u_{x}\wedge u_{y}$ in $L^1(\R^2)$, there exists a subsequence $n_{k_h}$ and a nonnegative function $g \in L^1(\R^2)$ such that $|u^{n_{k_h}}_{x}\wedge u^{n_{k_h}}_{y}| \leq g$ a.e. in $\R^2$. 
Thus, by the previous considerations and being $\varphi \in \R^3 + C_0^\infty(\R^2,\R^3)$ it holds that
$$|K(u^{n_{k_h}})-K(u)| |\varphi||u^{{n_{k_h}}}_{x}\wedge u^{{n_{k_h}}}_{y}| \leq 2C_K |\varphi|_\infty g,$$
and by the dominated convergence theorem we obtain that $I^{n_{k_h}}_1 \to 0$, contradicting $|I^{n_k}_1| > \epsilon$.
As far as concerns $I^n_2$, it suffices to observe that 
$$I^n_2 \leq C_K \|\varphi\|_\infty \int_{\R^{2}} |u^n_{x}\wedge u^n_{y} - u_{x}\wedge u_{y}| \to 0, \ \hbox{as} \ n \to \infty$$
because $u^n_{x}\wedge u^n_{y} \to u_{x}\wedge u_{y}$ in $L^1(\R^2)$. The proof is complete.
\qed

\begin{Remark}
\label{R:sphere2}
Let $\omega$ be the mapping introduced in Remark \ref{R:sphere}. Then, for every $\overline{p}\in\R^{3}$ and $r>0$, one has that $\Q(\overline{p}+r\omega)=-\int_{B_{r}(\overline{p})}K(p)~\!dp$, whereas if $r<0$ then $\Q(\overline{p}+r\omega)=\int_{B_{|r|}(\overline{p})}K(p)~\!dp$. For a proof of this fact see Remark 2.3 in \cite{CaMu11}.
\end{Remark}


Now we recall some useful results concerning the following volume-constrained minimization problems:
\begin{equation}
\label{eq:isovolumetric}
S_{K}(t):=\inf_{u\in M_{t}}\E(u)\quad\text{where}\quad\E(u):=\D(u)+\Q(u)~\!,
\end{equation}
$t\in\R$ is fixed, and $M_{t}$ is defined in (\ref{eq:Mt-def}). 
Unless differently specified, we always assume that $K\in C^{1}(\R^{3})$ satisfy $(K_{1})$ and $(K_{2})$. 

We point out that the mapping $t\mapsto S_{K}(t)$ is well defined on $\R$ and takes positive values for $t\ne0$, in view of (\ref{eq:isoperimetric}), (\ref{eq:QK-bound}), and (\ref{eq:VK-estimate}). It will be named the \emph{isovolumetric function}. 


\begin{Remark} 
\label{R:K=0}
For $t=0$ the class $M_{t}$ contains the constant functions. Since $0\le(1-\|Q_{K}\|_{\infty})\D(u)\le\E(u)$, we deduce that $S_{K}(0)=0$ and minimizers for $S_{K}(0)$ are exactly the constant functions. 
\end{Remark}
\begin{Remark} 
When $K=0$ we have $\E=\D$ and, by (\ref{eq:isoperimetric}), $S_{0}(t)=\inf\{\D(u)~|~u\in M_{t}\}=St^{2/3}$, for any fixed $t\in\R$. 
\end{Remark}

Now we state some properties of the isovolumetric function $S_{K}(t)$. 

\begin{Lemma}
\label{L:SK(t)}
For every $t\in\R$ the following facts hold:
\begin{itemize}[leftmargin=22pt]
\item[(i)] $S_{K}(-t)=S_{-K}(t)$;
\item[(ii)] $S_{K}(t)=S_{K(\cdot+p)}(t)$ for every $p\in\R^{3}$.
\item[(iii)] $S_{K}(t)=\inf\{\E(u)~|~u\in C^{\infty}_{c}(\R^{2},\R^{3}),~\V(u)=t\}$.
\end{itemize}
\end{Lemma}

\begin{Lemma}
\label{L:SK-estimates}
For every $t\in\R$ the following facts hold:
\begin{itemize}[leftmargin=20pt]
\item[(i)] For every $t\in\R$ one has that $(1-\|Q_{K}\|_{\infty})St^{2/3}\le S_{K}(t)\le S_{0}(t)=St^{2/3}$.
\item[(ii)] For every $t_{1},...,t_{k}\in\R$ one has that $S_{K}(t_{1})+...+S_{K}(t_{k})\ge S_{K}(t_{1}+...+t_{k})$.
\end{itemize}
\end{Lemma}

\begin{Remark}\label{R:star}
The value $S_{0}(t)$ is the infimum for the Dirichlet integral in the class $M_{t}$ of mappings in $\hat{H}^{1}$ parametrizing surfaces with volume $t$. We know that $S_{0}(t)$ is attained by a conformal parametrization of a round sphere of volume $t$ with arbitrary center (Lemma \ref{L:isoperimetric}). On the other hand, $S_{K}(t)$ is is the infimum value for the functional $\E=\D+\Q$ in the same class $M_{t}$, and $\Q$ has the meaning of $K$-weighted algebraic volume (see Remark \ref{R:sphere2}; see also \cite{CaMuARMA}, Sect.~2.3). 
\end{Remark}

The next result collects some properties about minimizing sequences for the isovolumetric problem defined by (\ref{eq:isovolumetric}). In particular we have a bound from above and from below on the Dirichlet norm, and we have that every minimizing sequence shadows another minimizing sequence consisting of approximating solutions for some prescribed mean curvature equation. 

\begin{Lemma}
\label{L:minimizing-PS}
Let $t\in\R$ be fixed. Then:
\begin{itemize}[leftmargin=22pt]
\item[(i)] 
$\D(u)\ge \frac{S_{K}(t)}{1+\|Q_{K}\|_{\infty}}$ for every $u\in M_{t}$.
\item[(ii)]
If $(u^{n})\subset M_{t}$ is a minimizing sequence for $S_{K}(t)$ then $\limsup\D(u^{n})\le \frac{St^{2/3}}{1-\|Q_{K}\|_{\infty}}$.
\item[(iii)]
For every minimizing sequence $(\tilde{u}^{n})\subset M_{t}$ for $S_{K}(t)$ there exists another minimizing sequence $(u^{n})\subset M_{t}$ such that $\|u^{n}-\tilde{u}^{n}\|\to 0$ and with the additional property that
\begin{equation*}
\Delta u^{n}-K(u^{n})u^{n}_{x}\wedge u^{n}_{y}+\lambda u^{n}_{x}\wedge u^{n}_{y}\to 0\quad\text{in }\hat{H}^{-1}\text{(= dual of $\hat{H}^{1}$)}
\end{equation*}
for some $\lambda\in\R$. 
\end{itemize}
\end{Lemma}


\begin{Definition}
Let $H\in C^{0}(\R^{3})$ be a given function. We call $U\in\hat{H}^{1}$ an $H$-bubble if it is a nonconstant solution to
\begin{equation}
\label{eq:Hsystem}
\Delta U=H(U)U_{x}\wedge U_{y}\quad\text{on}\quad\R^{2}
\end{equation} 
in the distributional sense. If $H$ is constant, an $H$-bubble will be named $H$-sphere. The system (\ref{eq:Hsystem}) is called $H$-system. 
\end{Definition}

A first useful property of $H$-bubbles, for a class of mappings $H$ of our interest, is the following:

\begin{Lemma}
\label{L:bubble-bdd}
Let $H(p)=K(p)-\lambda$ with $\lambda\in\R$ and $K\in C^{0}(\R^{3})$ satisfying $(K_{1})$. If $U\in\hat{H}^{1}$ is an $H$-bubble, then $U\in L^{\infty}$, and $\lambda\V(U)>0$. If, in addition, $K\in C^1(\R^3)$ then 
$U$ is of class $C^{2,\alpha}$ as a map on $\S^{2}$.  
\end{Lemma}

The next result is crucial and explains that Palais-Smale sequences for $\E$ constrained to $M_t$ admit a limit configuration made by bubbles. More precisely:

\begin{Lemma}[Decomposition Theorem]
\label{L:PS-decomposition}
Let $K\colon\R^{3}\to\R$ be a continuous function satisfying $(K_{1})$ and $(K_{2})$. If $(u^{n})\subset\hat{H}^{1}$ is a sequence  satisfying 
\begin{equation*}
\Delta u^{n}-K(u^{n})u^{n}_{x}\wedge u^{n}_{y}+\lambda u^{n}_{x}\wedge u^{n}_{y}\to 0\quad\text{in }\hat{H}^{-1},
\end{equation*}
 for some $\lambda\in\R$ and such that $c_{1}\le\|\nabla u^{n}\|_{2}\le c_{2}$ for some $0<c_{1}\le c_{2}<\infty$ and for every $n$, then there exist a subsequence of $(u^{n})$, still denoted $(u^{n})$, finitely many $(K-\lambda)$-bubbles $U^{i}$ ($i\in I$), finitely many $(-\lambda)$-spheres $U^{j}$ ($j\in J$) 
such that, as $n\to\infty$:
\begin{gather} \label{eqcPS}
\left\{\begin{array}{l}\D(u^{n})\to\textstyle\sum_{i\in I}\D(U^{i})+\sum_{j\in J}\D(U^{j})\\
\V(u^{n})\to\textstyle\sum_{i\in I}\V(U^{i})+\sum_{j\in J}\V(U^{j})\\
\Q(u^{n})\to\textstyle\sum_{i\in I}\Q(U^{i})
\end{array}\right.
\end{gather}
where $I$ or $J$ can be empty but not both. In particular, if $J=\varnothing$ then the subsequence $(u^{n})$ is bounded in $\hat{H}^{1}$. 
\end{Lemma}

\section{A constrained minimax result}
Let us denote by $\mathrm{Crit}_{\E}(t)$ the set of constrained critical points of the functional $\E$ over $M_{t}$, in the following sense:
$$
\mathrm{Crit}_{\E}(t)=\{u\in M_{t}~|~ \exists  \lambda \in \R \ \hbox{s.t.} \  \E^\prime(u)[\varphi] - \lambda \V^\prime(u)[\varphi] =0 \ \ \forall  \varphi \in \R^3+C_0^\infty(\R^2,\R^3) \}.
$$
For any $p \in \R^3$ and $t>0$ we set 
\begin{equation}\label{defwpt}
s_t:=\sqrt[3]{\frac{3t}{4\pi}}\quad\text{and}\quad
\omega_{p,t}:=s_t(-\omega + p),
\end{equation}
where $\omega$ the map defined in Remark \ref{R:sphere}.

The goal of this section is to prove the following result:
\begin{Proposition}\label{propcrucial}
Let $K \in C^1(\R^3)$ satisfying $(K_1)$--$(K_2)$ with \eqref{condkzero}, and $K>0$ on $\R^3$. 
Assume that 
\begin{equation}\label{assumptnex}
\exists t_0>0 \ \hbox{s.t.}~~ \mathrm{Crit}_{\E}(t)=\varnothing\ \ \forall  t \in (0,t_0].
\end{equation}
Then there exists $R>0$ such that for every $t\in(0,t_{0})$
$$
S_0(t) < \sup_{p \in \partial B_R} \E(\omega_{p,t}) < \inf_{\phi \in \Phi} \sup_{p \in \overline{B_R}} \E(\phi(p)) < 2^{1/3} S_0(t),
$$
where $S_0(t)=St^{2/3}$,  $\Phi:=\{\phi \in C^0(\overline{B_R},M_t); \ \phi|_{\partial B_R}(p)=\omega_{p,t}\}$, $\omega_{p,t}$ is the function defined in \eqref{defwpt}.
\end{Proposition}

In order to prove Proposition \ref{propcrucial} we need to introduce a new tool and some preliminary results. Let us fix $t>0$ and denote by $\B_t\colon \hat{H}^{1} \to \R^3$ the vector-valued map defined by
$$\B_t(u):=\frac{1}{8\pi s_t^2}\int_{\R^2} \Pi(u) |\nabla u|^2 ,$$
where $\Pi$ is the minimal distance projection of $\R^3$ onto the closed unit ball, namely 
$$\Pi(p):=\begin{cases}p & \hbox{if} \ |p|<1 \\ \frac{p}{|p|} & \hbox{if} \ |p|\geq 1.\end{cases}.$$
Since $\Pi\circ u$ is bounded for any $u \in \hat{H}^{1}$, the mapping $\B_t$ is well defined and continuous on $\hat{H}^1$, in particular $\B_t$ is continuous as mapping from $M_{t}$ to $\R^3$. We also point out that $\B_t$ is conformally invariant.
 
 \begin{Proposition}\label{Propfond2}
 Let $K \in C^1(\R^3)$ satisfying $(K_1)$--$(K_2)$ and assume that $K>0$ in $\R^3$ and \eqref{assumptnex}. Then
 $$ \forall t \in(0,t_0)  \ \exists r_{t,K}>0~~\hbox{s.t.} \ \inf_{\scriptstyle{u \in M_t}\atop\scriptstyle{|\B_t(u)|\leq r_{t,K}}} \E(u)> S_0(t). $$
 \end{Proposition}
\Proof
We argue by contradiction. Assume the thesis is false, then there exist $t \in (0,t_0)$ and a sequence $(u^n) \subset M_t$ such that $\E(u^n) \to S_0(t)$ and $\B_t(u^n)\to 0$. Observe that, thanks to Theorem \ref{T:nonexistence}, we can assume without loss of generality that $S_K(t)=S_0(t)$. Hence $(u^n)$ is a minimizing sequence for $S_K(t)$ and, by Lemma \ref{L:minimizing-PS}, there exists another minimizing sequence $(\tilde u^n) \subset M_t$ such that $\|u^n - \tilde u^n\| \to 0$ and
$$
\Delta u^{n}-K(u^{n})u^{n}_{x}\wedge u^{n}_{y}+\lambda u^{n}_{x}\wedge u^{n}_{y}\to 0\quad\text{in }\hat{H}^{-1}
$$
for some $\lambda\in\R$.
Now, being $(\D(u^n))$ bounded, by Lemma \ref{L:PS-decomposition}, we get that, up to a subsequence (still denoted $(u^n)$), there exist finitely many $(K-\lambda)$-bubbles $U^{i}$ ($i\in I$), finitely many $(-\lambda)$-spheres $U^{j}$ ($j\in J$) for which \eqref{eqcPS} holds, and $I$ or $J$ can be empty but not both. Since we are assuming \eqref{assumptnex}, it results that $I=\varnothing$ and thus $J\neq \varnothing$. Now we prove that $J$ is a singleton. Assume, by contradiction, that $J$ is not a singleton, in particular, being $J$ finite and denoting by $|J|$ its cardinality, we have $|J|\geq 2$. We set $t_j:=\V(U_j)$ for $j \in J$. By Lemma \ref{L:bubble-bdd}  one has that $t_j \lambda>0$ for any $j \in J$. Hence, from \eqref{eqcPS} we get that $\sum_{j \in J} t_j=t$, and  $t_j >0$ for any $j \in J$. We observe that for any $j \in J$ being $U_j$ a $(-\lambda)$-sphere, there exists a positive integer $k_j$ such that
\begin{gather} \label{eqalgebraicvolum}
\begin{array}{l}
4\pi k_j \lambda^2 = \D(U_j),\ \frac{4}{3} \pi k_j \lambda^3=t_j.
\end{array}
\end{gather}
From \eqref{eqalgebraicvolum} we deduce that $\D(U_j)=Sk_j^{1/3}t_j^{2/3}$.
Moreover, thanks to \eqref{eqcPS}, being $S_K(t)=S_0(t)$ we have $S_0(t)=\sum_{j \in J} \D(U_j)$, and being $k_j \in \N^+$ it holds
$$
\Big(\sum_{j \in J} t_j\Big)^{2/3}=\sum_{j \in J} k_j^{1/3} t_j^{2/3} \geq \sum_{j \in J} t_j^{2/3}.
$$
On the other hand, being $t_j>0$ for all $j \in J$ and $|J|\geq 2$, by a well known elementary inequality, it also holds
\begin{equation*}
 \Big(\sum_{j \in J} t_j\Big)^{2/3}< \sum_{j \in J} t_j^{2/3},
\end{equation*}
which gives a contradiction. Now, being $J$ a singleton, by Theorem 0.1 of \cite{CaJFA} and thanks to \eqref{assumptnex}, there exists a sequence $(g_n)$ of conformal transformations of $\R^2 \cup \{\infty\}$ into itself such that setting
$$v_n: = \tilde u^n \circ g_n \quad\hbox{and}\quad p_{n}:= \frac{1}{4\pi}\int_{\R^2} \mu^2 v_n$$
for a subsequence of $(v_n)$ (still denoted by $(v_n)$), one has $|p_n|\to \infty$ and $v_n - p_n \to U_j  \ \hbox{weakly in} \ \hat{H}^1.$ In particular, being $\D(\tilde u^n) \to \D(U_j)$ it holds that $\nabla v_n \to \nabla U_j \ \hbox{in} \ L^2.$
Recalling that $\B_t$ is conformally invariant we get that
$$
\B_t(\tilde u^n)=\B_t(v_n)=\frac{1}{8\pi s_t^2} \int_{\R^2} \Pi(v_n) |\nabla U_j|^2 + o(1).
$$
Since $|p_n|\to\infty$ and $\int_{\S^2} U_j = 0$ we also have that $v_n - p_n \to U_j$ strongly in $\hat{H}^1$. In particular $|v_n|\to\infty$ a.e. in $\R^2$. Being $p_n/|p_n|$ bounded, up to a subsequence, we have $p_n/|p_n| \to p \in \S^2$ and it follows that $\Pi(v_n) \to p$. Hence we obtain that
$$|\B_t(\tilde u^n)| = \frac{1}{8\pi s_t^2} \int_{\R^2} |\nabla U_j|^2 + o(1) \geq c>0,$$
and being $\B_t$ continuous and $\|\tilde u^n-u^n\|\to 0$, this contradicts $\B_t(u^n)\to 0$.
 \qed
 
\begin{Lemma}\label{lemmasemplic1}
Let $K \in C^1(\R^3)$ satisfying $(K_1)$, $(K_2)$ and \eqref{assumptnex}. Let $t \in (0,t_0)$ and $p \in \R^3$. As $|p|\to\infty$ it holds that
$$\E(\omega_{p,t})=S_0(t)+o(1),$$
$$ \B_t(\omega_{p,t})=\frac{p}{|p|}+o(1),$$
where $\omega_{p,t}$ is the function defined in \eqref{defwpt}.
\end{Lemma}
\Proof
The first relation follows from the fact that $\D(\omega_{p,t})=S_0(t)$ and $\Q(\omega_{p,t})=\int_{B_{s_t}(s_tp)}K(q) \ dq$ (see Remark \ref{R:sphere2}). Thanks to assumption $(K_1)$ we get that
$$|\Q(\omega_{p,t})| \leq \int_{B_{s_t}(s_tp)}|K(q)| \ dq \leq k_0 \int_{B_{s_t}(s_tp)} \frac{dq}{|q|}.$$ Recalling that $1/|q|$ is harmonic in $\R^3$ outside the origin we have
$$ \int_{B_{s_t}(s_tp)} \frac{dq}{|q|} = \frac{4\pi s_t^2}{3|p|} \to 0 \ \hbox{as} \ |p|\to\infty.$$
The first relation is then proved. 
Concerning the second relation, we observe that $|\omega_{p,t}|\geq 1$ on $\R^2$ for $|p|$ large enough. This implies that 
$\B_t(\omega_{p,t})=s_t^2 \B_t(-\omega+p)$ and
$$s_t^2 \B_t(-\omega+p)-\frac{p}{|p|}=\frac{1}{8\pi} \int_{\R^2} \left(\frac{\omega+p}{|\omega+p|} - \frac{p}{|p|}\right) |\nabla \omega|^2\to 0 \ \hbox{as} \ |p|\to\infty,$$
by dominated convergence theorem. The proof is then complete.
\qed

We have now all the tools to prove Proposition \ref{propcrucial}.\\

\Proof
Let $t \in (0,t_0)$, let $r_{t,K}>0$ be given by Proposition \ref{Propfond2} and let $\epsilon \in (0,1)$ be such that $\epsilon < \inf_{u \in M_t, |\B_t(u)|\leq r_{t,K}} \E(u)- S_0(t)$. According to Lemma \ref{lemmasemplic1} there exists a sufficiently large number $R>1$  such that for all $p \in \R^3$ with $|p|=R$ one has that
\begin{equation}\label{eq1dimpropcrucial}
\E(\omega_{p,t}) < S_0(t) + \epsilon, \ \ |\B_t(\omega_{p,t})|>1-\epsilon, \ \ \frac{p}{|p|} \cdot \B_t(\omega_{p,t}) > 1 -\epsilon.
\end{equation}
Let $\Phi$ be as in the statement of Proposition \ref{propcrucial}. Being $K>0$, and thanks to \eqref{eq1dimpropcrucial}, it follows that:
$$S_0(t) < \sup_{p \in \partial B_R} \E(\omega_{p,t}) < S_0(t) + \epsilon.$$
Let us set $c:=\inf_{\phi \in \Phi} \sup_{p \in \overline{B_R}} \E(\phi(p))$.
We want to prove that $ c\geq S_0(t) + \epsilon$. To this goal, assume by contradiction that there exists a map $\phi \in \Phi$ such that
$$\sup_{p \in \overline{B_R}} \E(\phi(p)) < S_0(t) + \epsilon.$$
Hence by Proposition \ref{Propfond2} we have that 
\begin{equation}\label{eqconsPropfond2}
|\B_t(\phi(p))| > r_{t,K} \ \hbox{for all} \ p \in \overline{B_R}.
\end{equation}
Now consider the map $g:\overline{B_R} \to \R^3$ defined by
$$g(p):=\B_t(\phi(p)),$$
and fix a point $p_0 \in \R^3$ with $0<|p_0|< \min\{r_{t,K},1-\epsilon\}$. We claim that the topological degree $\mathrm{deg}(g,\overline{B_R},p_0)$ is well defined and $\mathrm{deg}(g,\overline{B_R},p_0)=1.$
To this purpose, consider the homotopy $h\colon[0,1] \times \overline{B_R} \to \R^3$ defined by
$$h(s,p):=sp+(1-s)\B_t(\phi(p)).$$
Assume by contradiction that $h$ is not admissible, then, there exist $\bar s \in [0,1]$ and $\bar p \in \partial B_R$ such that $h(\bar s, \bar p)=p_0$, hence by definition of $h$ and thanks to \eqref{eq1dimpropcrucial} we deduce that
$$|p_0|=|\bar s \bar p + (1-\bar s) \B_t(\omega_{\bar p, t})| \geq \bar s R + (1-\bar s) (1-\epsilon)=1-\epsilon + \bar s (R-1+\epsilon)\geq 1-\epsilon,$$
which gives a contradiction. Hence $h$ is an admissible homotopy between $g$ and the identity map of $\overline B_R$, and by well known properties of the topological degree it holds that $\mathrm{deg}(g,\overline{B_R},p_0)=1$.
Now, being $\mathrm{deg}(g,\overline{B_R},p_0)\neq 0$ in particular we deduce that the equation $g(p)=p_0$ has at least a solution $p \in B_R$. Hence $|\B_t(\phi(p))|=|p_0|$ for some $p \in B_R$ but, being $|p_0|<  r_{t,K}$ it follows that $|\B_t(\phi(p))| <  r_{t,K}$, 
contradicting \eqref{eqconsPropfond2} and hence $ c\geq S_0(t) + \epsilon$.

In order to conclude the proof, it remains to check that $c<2^{1/3} S_0(t)$. To this goal, let us consider the map $p \mapsto \omega_{p,t}$. It is clear that $\omega_{p,t} \in \Phi$. It is known that $\D(\omega_{p,t})=S_0(t)$, hence in order to complete the proof we need to estimate $\Q(\omega_{p,t})$. By Remark \ref{R:sphere2}, we know that
 $\Q(\omega_{p,t})=\int_{B_{s_t}(s_tp)}K(q) \ dq$. Thanks to assumption $(K_1)$ and being $K>0$, it holds that
 $$\int_{B_{s_t}(s_tp)}K(q) \ dq \leq k_0 \int_{B_{s_t}(s_tp)} \frac{1}{|q|}\ dq~\!.$$
 By a suitable change of variable and elementary computations we get that
 $$k_0\int_{B_{s_t}(s_tp)}K(q) \ dq =k_0 s_t^2 \int_{B_1(0)} \frac{1}{|q-p|} \ dq.$$ 
 Let us set consider the function $I: \R^3 \to \R$ defined by $p \mapsto  \int_{B_1(0)} \frac{1}{|q-p|}$.
We observe that $I(p)$ can be explicitly computed, more precisely we have that
\begin{equation}\label{eq1formexpl}
I(p) = \begin{cases} \frac{2\pi}{3}(3-|p|^2) & \hbox{if} \ |p|\leq 1,\\
\frac{4\pi}{3|p|} &  \hbox{if} \ |p|> 1. \end{cases}
\end{equation}
In fact if $|p|>1$ the integrand function $ q \mapsto \frac{1}{|q-p|}$ is harmonic in $B_1(0)$ and thus by the mean value property we get that $I(p)= \frac{4\pi}{3|p|}$.

 The case $|p|\leq 1$ is more delicate: first observe that by dominated convergence theorem we get that $I:\R^3 \to \R$ is continuous at any $p \in \R^3$, in particular if $p \in \partial B_1(0)$, from the previous case, we deduce that $I(p)=\frac{4}{3} \pi$.
Now let us consider the vector field $E\colon\R^3 \to \R^3$ defined by $$ E(p):=\int_{B_1(0)}\frac{p-q}{|p-q|^3} \ dq.$$
Observe that, by definition and by the dominated convergence theorem, $I$ is differentiable and
\begin{equation}\label{gradenergypot}
E(p)=-\nabla I (p).
\end{equation}
We also note that $E$ is of the form $E(p)=g(|p|) \frac{p}{|p|}$ when $p\neq 0$, for some function $g\colon\R^+ \to \R$. In fact, fixing $p\neq 0$ and making a change of variable in the integral defining $E$ by any orthogonal matrix $T \in O_3$ such that $T(p)=p$ we get that $E(p)=T(E(p))$, and thus the fact follows from the arbitrariness of $T$. Moreover, since $T(0)=0$ for any $T \in O_3$ it holds that $E(0)=0$. At the end, by a suitable application of the Stokes Theorem, we obtain that
 $$ E(p) = \begin{cases} \frac{4}{3}\pi p & \hbox{if} \ |p|\leq 1,\\  \frac{4}{3}\pi {p}/{|p|^3} & \hbox{if} \ |p|> 1. \end{cases}$$
Thanks to \eqref{gradenergypot} and the previous characterization, by fixing a point $p_0 \in \partial B_1$, and for any path $\gamma$ joining $p$ and $p_0$, we get that
$$I(p) - I(p_0) = - \int_{\gamma} E \cdot d\gamma =-   \frac{4}{3}\pi \int_1^{|p|} r \ dr=-  \frac{2}{3}\pi(|p|^2-1).$$
Hence, since $I(p_0)= \frac{4}{3}\pi$ we have
$$I(p)=   \frac{4}{3}\pi  -  \frac{2}{3}\pi(|p|^2-1) = \frac{2}{3}\pi(3-|p|^2)= \frac{2}{3}\pi (3-|p|^2).$$
 Hence, thanks to \eqref{eq1formexpl} we deduce that $$\sup_{p \in \R^3} I(p)= 2\pi$$ and
$$k_0 s_t^2 \int_{B_1(0)} \frac{1}{|q-p|} \ dq \leq k_0 s_t^2 2\pi.$$ 
Now observe that $$\E(\omega_{p,t}) \leq S_0(t) + 2 \pi k_0 \left(\frac{3}{4\pi}\right)^{2/3} t^{2/3}= S t^{2/3} + 2k_0 \left(\frac{9}{16}\pi\right)^{1/3} t^{2/3}.$$ 
Now,  thanks to the assumption \eqref{condkzero}, by elementary computations it is easy to verify that $$S t^{2/3} + 2k_0 \left(\frac{9}{16}\pi\right)^{1/3} t^{2/3} < 2^{1/3} S t^{2/3},$$ which implies that $\E(\omega_{p,t}) <  2^{1/3}S t^{2/3}$, and in particular $c< 2^{1/3} S t^{2/3}$ which is the desired relation. The proof is complete.
\qed

\section{Constrained Palais-Smale sequences for $\E$ at the minimax level $c$}
In this section we prove that there exists a Palais-Smale sequence constrained to the smooth manifold $M_t$ at a suitable minimax level. Let $t>0$ and $R>0$ be fixed, we define
\begin{equation}\label{cminimaxlevel}
c:=\inf_{\phi \in \Phi} \sup_{p \in \overline{B_R}} \E(\phi(p)),
\end{equation}
where $\Phi:=\{f \in C^0(\overline{B_R},M_t); \ f|_{\partial B_R}(p)=\omega_{p,t}\}$, $\omega_{p,t}$ is the function defined in \eqref{defwpt}. Moreover we define
 \begin{equation}\label{czero}
 c_0:= \sup_{p \in \partial B_R} \E(\omega_{p,t}).
\end{equation}
We begin with a preliminary result:

\begin{Lemma} \label{lemmadensity}
Let $t \in \R\setminus\{0\}$ and let $u \in M_t$. It holds that $T_u M_t \cap (\R^3 + C_0^\infty(\R^2,\R^3))$ is dense in $T_u M_t$.
\end{Lemma}
\Proof
Let us fix $t \in \R\setminus\{0\}$. By Lemma \ref{L:volume} we know that $M_t \subset \hat{H}^1$ is a smooth manifold of codimension one. Let us fix $u \in M_t$. Then we can write $$\hat{H}^1=T_uM_t \oplus \langle h \rangle,$$
where $h \in \hat{H}^1$ is the Riesz rapresentative of $\frac{\V^\prime(u)}{\|\V^\prime(u)\|^2}$ (see Section 6.1 of \cite{AmbMalc}).
We observe that since $\R^3 + C_0^\infty(\R^2,\R^3)$ dense in $\hat{H}^1$ (see Lemma \ref{L:density}) then there exists $v \in (\R^3 + C_0^\infty(\R^2,\R^3)) \setminus T_uM_t$. Hence we can also write
$$\hat{H}^1=T_uM_t \oplus \langle v \rangle.$$
Now let us fix $w \in T_uM_t$, then by the density of $\R^3 + C_0^\infty(\R^2,\R^3)$ in $\hat{H}^1$ there exists a sequence $(w_n) \subset \R^3 + C_0^\infty(\R^2,\R^3)$ such that $w_n \to w$ in $\hat{H}^1$. 
Let us set $$\tilde w_n:= w_n - \frac{\V^\prime(u)[w_n]}{\V^\prime(u)[v]} v.$$
By construction, $\tilde w_n \in \R^3 + C_0^\infty(\R^2,\R^3)$ and $\V^\prime(u)[\tilde w_n]=0$, i.e. $\tilde w_n \in (\R^3 + C_0^\infty(\R^2,\R^3)) \cap T_uM_t$. Moreover
$$\| \tilde w_n - w\| \leq \|w_n -w\| + \left|\frac{\V^\prime(u)[w_n]}{\V^\prime(u)[v]} \right|\|v\|,$$
and the right-hand side goes to zero as $n\to\infty$ because $w_n \to w$ in $\hat{H}^1$ and $w \in T_uM_t$. The proof is complete.
\qed

\begin{Proposition}\label{propconstrminmaxPS}
 Let $t \in \R^+$ and $R>0$ be fixed and let $c$, $c_0$ be the numbers defined, respectively, in \eqref{cminimaxlevel}, \eqref{czero}. If $c>c_0$ then for any sufficiently small $\epsilon >0$ and for each $f \in \Phi$ such that 
\begin{equation}\label{eq1propconstrPS}
\sup_{p \in \overline{B_R}} \E(f(p)) \leq c + \epsilon
\end{equation}
there exists $u \in M_t$ such that
\begin{gather*}
c-\epsilon \le\E(u)\le\sup_{p \in \overline{B_R}} \E(f(p))~\!,\\
\|u-f(p)\|\le\epsilon^{1/2}~~\forall p\in\overline{B_{R}}~\!,\\
|\E^\prime(u)[\varphi]| \leq 2 \epsilon^{1/2}~~\forall \varphi \in T_u M_t \cap (\R^3 + C_0^\infty(\R^2,\R^3)) \ \hbox{with} \ \|\varphi\|=1~\!.
\end{gather*}
\end{Proposition}
\Proof
Let $\epsilon$ be such that $0<\epsilon < c-c_0$. Moreover assume that $\epsilon$ satisfies
\begin{equation}\label{restrictiononeps}
\epsilon^2 \left(\frac{1}{3t} + \frac{2}{9}2^{7/3} \epsilon^2\right)<1.
\end{equation}
 A further restriction on the smallness of $\epsilon$ will be specified in the sequel of the proof. Let $f \in \Phi$ satisfy \eqref{eq1propconstrPS} and define the function $F:\Phi \to \R$ by setting
$$F(g):=\sup_{p \in \overline{B_R}} \E(g(p)).$$
In particular observe that $c=\inf_{\Phi} F > c_0$. Thanks to Ekeland's variational principle (see, e.g., \cite{Eke79}) there exists $h \in \Phi$ such that
\begin{gather}
\label{eq1bispropconstrPS}
F(h) \leq F(f) \leq  c + \epsilon~\!,
\\
\nonumber
d(h,f):=\sup_{p \in \overline{B_R}} \|h(p)-f(p) \| \le \epsilon^{1/2}~\!,
\\
\nonumber
F(g)>F(h) - \epsilon^{1/2} d(h,g)~~\forall g \in \Phi \ \hbox{with} \ g\neq h~\!.
\end{gather}
In order to reach the conclusion, it suffices to show that for some $p \in \overline{B_R}$ it holds that
\begin{gather}
\nonumber
c-\epsilon \leq \E(h(p))~\!,\\
\label{eq:Eprime}
|\E^\prime(h(p))[\varphi]| \leq 2 \epsilon^{1/2}~~\forall \varphi \in T_{h(p)} M_t \cap (\R^3 + C_0^\infty(\R^2,\R^3)) \ \hbox{with} \ \|\varphi\|=1~\!.
\end{gather}
Notice that \eqref{eq:Eprime}  is equivalent to 
$$\E^\prime(h(p))[\varphi] \geq - 2 \epsilon^{1/2}~~\forall\varphi \in T_u M_t \cap (\R^3 + C_0^\infty(\R^2,\R^3)) \ \hbox{with} \ \|\varphi\|=1~\!.
$$
By contradiction, if this does not happen, then, setting
$$P:=\{p \in \overline{B_R}~|~ c-\epsilon \leq \E(h(p)) \},$$
for each $p \in P$ there exists $\delta_p>0$, $\varphi_p \in  T_{h(p)} M_t \cap (\R^3 + C_0^\infty(\R^2,\R^3))$ with $\|\varphi_p\|=1$ and an open ball $B_p$ centered at $p$ such that for $q \in B_p$ and $u \in \hat{H}^1$ with $\|u\| \leq \delta_p$, we have
\begin{equation}\label{eq2propconstrPS}
 \E^\prime(h(q)+u)[\varphi_p] < - 2 \epsilon^{1/2}.
\end{equation}
We recall that by Lemma \ref{techlemmacont} for any fixed $\varphi \in \R^3 + C_0^\infty(\R^2,\R^3)$ the map $u \mapsto \E^\prime(u)[\varphi]$ from $\hat{H}^1$ to $\R$ is continuous.
Moreover, since $\V^\prime(h(p))[\varphi_p]=0$, taking a possibily smaller ball $B_p$ and a smaller constant $\delta_p$, if necessary, we can also assume that 
\begin{equation}\label{eq2bispropconstrPS}
|\V^\prime(h(q)+u) [\varphi_p] | \leq \epsilon^2~~\forall q\in B_{p}~\!,~~\forall u \in \hat{H}^1\text{ with $\|u\| \leq \delta_p$~\!.}
\end{equation}
By the continuity of $\D$, assumption $(K_1)$ and \eqref{eq1bispropconstrPS}, taking a smaller $\delta_p$, we can also assume that 
\begin{equation}\label{eq2quaterpropconstrPS}
\D(h(p) + u) \leq C~~\text{for $\|u\|\leq \delta_p$~\!,}
\end{equation}
where $C$ is some positive constant depending only on $k_0$ and $c$.
In fact, by assumption $(K_1)$ we have  
$$
\D(h(p)) \leq \frac{\E(h(p))}{1-k_0/2} \leq \frac{F(h)}{1-k_0/2} \leq \frac{c + \epsilon}{1-k_0/2} < C,
$$
for some positive constant $C$ depending only on $k_0$ and $c$ and then by continuity of $\D$ we get \eqref{eq2quaterpropconstrPS}.
Since $P$ is compact there exists a finite subcovering $B_{p_1},\ldots,B_{p_k}$ of $P$ and we define $\psi_j\colon P\to [0,1]$ by 
\begin{equation*}
\psi_j(p)=\left\{\begin{array}{ll}
\displaystyle \frac{\mathrm{dist}(p, \complement B_{p_j})}{\sum_{i=1}^k \mathrm{dist}(p, \complement B_{p_i})}&  \hbox{if~} \ p \in \bigcup_{i=1}^k B_{p_i},\\[8pt]
0 &\hbox{if~} \ p \in P \setminus \bigcup_{i=1}^k B_{p_i}.
\end{array}\right.
\end{equation*}
Furthermore let $\delta:=\min\{\frac{1}{2},\frac{t}{2},\delta_{p_1}, \ldots, \delta_{p_k}\}$, let $\psi\colon\overline{B_R} \to [0,1]$ be a continuous function such that
\begin{equation*}
\psi(p)=\left\{\begin{array}{ll}
1&\hbox{if~} \ c\leq \E( h(p)),\vspace{4pt}\\
0&\hbox{if~} \  \E(h(p)) \leq c - \epsilon,
\end{array}\right.
\end{equation*}
and let $\tau\colon\overline{B_R} \to \R$ and $g\colon\overline{B_R}\to M_t$ be defined by
$$
\tau(p):=\sqrt[3]{\frac{t}{\V\left(h(p) + \delta \psi(p) \sum_{j=1}^k \psi_j(p)\varphi_{p_j}\right)}}~\!,\quad
g(p):=\tau(p) \left(h(p) + \delta \psi(p) \sum_{j=1}^k \psi_j(p)\varphi_{p_j} \right)~\!.
$$
It holds that $g \in \Phi$. In fact, since $0<\epsilon<c-c_0$, when $p \in \partial B_R$ we have
$$\E(h(p)) = \E(s_t(-\omega + p)) \leq c_0 < c-\epsilon,$$
and hence $\psi(p)=0$ which means that $g(p)=h(p)=s_t(-\omega+ p)$.
We observe that for $p \in P$ the following inequality holds: 
\begin{equation}\label{eq2trispropconstrPS}
1- \frac{1}{3t}\delta \psi(p) \epsilon^2 - \frac{2}{9} 2^{7/3}\delta^2 \psi^2(p)\epsilon^4 \leq\tau(p) \leq 1 + \frac{1}{3t}\delta \psi(p) \epsilon^2 + \frac{2}{9} 2^{7/3}\delta^2 \psi^2(p)\epsilon^4.
\end{equation}
In fact, by the mean value theorem we have
\begin{equation*}
\begin{split}
\V\left(h(p) + \delta  \psi(p) \sum_{j=1}^k \psi_j(p)\varphi_{p_j}\right)&=\V(h(p)) + \V^\prime \left(h(p) +  \sigma \delta  \psi(p) \sum_{j=1}^k \psi_j(p)\varphi_{p_j}\right) [ \delta \psi(p) \sum_{j=1}^k \psi_j(p)  \varphi_{p_j}]\\
&=t + \delta \psi(p) \sum_{j=1}^k \psi_j(p)  \V^\prime \left(h(p) +  \sigma \delta  \psi(p) \sum_{j=1}^k \psi_j(p)\varphi_{p_j}\right) [\varphi_{p_j}]
\end{split}\end{equation*}
for some $\sigma \in (0,1)$. Now, thanks to \eqref{eq2bispropconstrPS} and the definition of the functions $\psi_j$, we see that
$$\left|\sum_{j=1}^k \psi_j(p)  \V^\prime \left(h(p) +  \sigma \delta  \psi(p) \sum_{j=1}^k \psi_j(p)\varphi_{p_j}\right) [\varphi_{p_j}]\right| \leq \epsilon^2.$$
In particular we observe that this estimate is uniform with respect to $p \in P$. Hence we deduce that
$$\tau(p) =\sqrt[3]{\frac{t}{t + \delta \psi(p) O(\epsilon^2)}}$$
with $|O(\epsilon^2)|\leq \epsilon^2$ and the desired inequality follows by elementary considerations. More precisely, by the Taylor expansion of the function $s\mapsto\frac{1}{(1+s)^{1/3}}$ we have
 $$\tau(p) =1 - \frac{1}{3t} \delta \psi(p) O(\epsilon^2) + \int_0^1(1-s)\frac{4}{9}\left(1+ s \left(\frac{\delta}{t} \psi(p) O(\epsilon^2)\right)\right)^{-7/3} \left(\frac{\delta}{t} \psi(p) O(\epsilon^2)\right)^2~\! ds.$$
Thanks to the choice of $\delta$ and being $|O(\epsilon^2)|\leq \epsilon^2$  we have $|\frac{\delta}{t} \psi(p) O(\epsilon^2)| \leq \frac{1}{2}\epsilon^2 \leq \frac{1}{2}$. Hence, for any $s\in [0,1]$ we have $\left(1+ s\left(\frac{\delta}{t} \psi(p) O(\epsilon^2)\right)\right)^{-7/3} \leq 2^{7/3}$ and we get that
$$\left|\int_0^1(1-s)\frac{4}{9}\left(1+ s \left(\frac{\delta}{t} \psi(p) O(\epsilon^2)\right)\right)^{-7/3} \left(\frac{\delta}{t} \psi(p) O(\epsilon^2)\right)^2 \!~ ds\right| \leq \frac{2}{9} 2^{7/3} \frac{\delta^2}{t^2} \psi^2(p) \epsilon^4.$$
Hence the estimate \eqref{eq2trispropconstrPS} follows immediately.
Now, setting $\eta(p):= \psi(p) \sum_{j=1}^k \psi_j(p)\varphi_{p_j} $, we write
\begin{eqnarray*}
\E(g(p))- \E(h(p)) & =&\underbrace{\tau(p)^2 \left(\E (h(p) + \delta \eta(p)) - \E(h(p))\right)}_{I_1} + \underbrace{\tau(p)^2 \E (h(p)) - \E(h(p))}_{I_2} \\
&&+ \underbrace{\Q(\tau(p)( h(p)+ \delta \eta(p))) - \Q(h(p) + \delta \eta(p))}_{I_3}\\
&& + \underbrace{\Q(h(p) + \delta \eta(p)) - \tau(p)^2 \Q(h(p) + \delta \eta(p))}_{I_4}  .
\end{eqnarray*}
We begin with the term $I_1$. Recalling that for any fixed $\varphi \in \R^3 + C_0^\infty(\R^2,\R^3)$ the functional $\E$ is  differentiable along $\varphi$, by the mean value theorem, for any fixed $p \in P$ there exists $\xi \in (0,1)$ such that
\begin{equation}\label{eq3propconstrPS}
\begin{split}
\E( h(p) +  \delta \eta(p)) - \E(h(p))&=\displaystyle \E^\prime\left( h(p)+ \xi \delta \psi(p) \sum_{j=1}^k \psi_j(p)\varphi_{p_j}\right)[ \delta \psi(p) \sum_{j=1}^k \psi_j(p)\varphi_{p_j}]\\
&= \displaystyle  \delta \psi(p) \sum_{j=1}^k \psi_j(p) \E^\prime\left(h(p)+ \xi \delta \psi(p) \sum_{j=1}^k \psi_j(p)\varphi_{p_j}\right)[\varphi_{p_j}]\\
&= \displaystyle  \delta \psi(p) \sum_{j=1}^k \psi_j(p) \E^\prime\left(h(p)+ \xi\tau (p)\delta \psi(p) \sum_{j=1}^k \psi_j(p)\varphi_{p_j}\right)[\varphi_{p_j}]~\!.m
\end{split}
\end{equation}
Now, from \eqref{eq2propconstrPS}, \eqref{eq2trispropconstrPS} and \eqref{eq3propconstrPS} we get that for $p \in P$
$$
I_{1}\leq - 2\tau(p)^2 \delta \psi(p) \epsilon^{1/2}\leq - 2  \left(1-\frac{1}{3t} \delta \psi(p) \epsilon^2 - \frac{1}{9t^2} 2^{7/3} \delta^2 \psi^2(p) \epsilon^4\right)^2 \delta \psi(p) \epsilon^{1/2}~\!.
$$
Regarding the term $I_2$, thanks to \eqref{restrictiononeps}, \eqref{eq1bispropconstrPS} and  \eqref{eq2trispropconstrPS} we have
$$|I_2|= |\tau^2(p)-1| \E(h(p)) \leq 3 (c+\epsilon)  \left(\frac{1}{3t} \delta \psi(p) \epsilon^2 + \frac{2}{9t^2} 2^{7/3} \delta^2 \psi^2(p) \epsilon^4\right).$$
For $I_3$, thanks to Lemma \ref{L:K-volume}, we have that
\begin{equation*}
\begin{split}
&\Q(\tau(p)( h(p)+\delta \eta(p))) -\Q(h(p)+\delta \eta(p))\\
&\qquad=\int_1^{\tau(p)} s^2 \int_{\R^2}K(s  (h(p)+\delta \eta(p)))h(p)+\delta \eta(p) \cdot ( h(p)+\delta \eta(p))_x \wedge ( h(p)+\delta \eta(p))_y. 
\end{split}\end{equation*}
Now, by assumption $(K_1)$ and thanks to \eqref{restrictiononeps}, \eqref{eq2quaterpropconstrPS}, \eqref{eq2trispropconstrPS} we get that
\begin{equation*}
\begin{split}
& \left| \int_1^{\tau(p)} s^2 \int_{\R^2}K(s  (h(p)+\delta \eta(p))) h(p)+\delta \eta(p) \cdot ( h(p)+\delta \eta(p))_x \wedge ( h(p)+\delta \eta(p))_y\right| \\
&\qquad= \left| \int_1^{\tau(p)} s \int_{\R^2}K(s  (h(p)+\delta \eta(p))) s h(p)+\delta \eta(p) \cdot ( h(p)+\delta \eta(p))_x \wedge ( h(p)+\delta \eta(p))_y\right| \\
 &\qquad\leq \int_{\min\{1,\tau(p)\}}^{\max\{1,\tau(p)\}}  s \int_{\R^2} \left|K(s  (h(p)+\delta \eta(p))) s (h(p)+\delta \eta(p)) \right| \left|( h(p)+\delta \eta(p))_x \wedge ( h(p)+\delta \eta(p))_y\right| \\
 &\qquad\leq k _0 \int_{\R^2} \left|( h(p)+\delta \eta(p))_x \wedge ( h(p)+\delta \eta(p))_y\right| \left| \int_1^{\tau(p)}  s  \ ds\right|\\
  &\qquad\leq k _0 \D(h(p)+\delta \eta(p)) \frac{|\tau(p)^2 -1|}{2}\\
  &\qquad\leq\frac{3}{2} k_0 C  \left(\frac{1}{3t} \delta \psi(p) \epsilon^2 + \frac{2}{9t^2} 2^{7/3} \delta^2 \psi^2(p) \epsilon^4\right)~\!.
\end{split}\end{equation*}
As far as concerns $I_4$, as before, using assumption $(K_1)$ we get that
\begin{eqnarray*}
|I_4|&=&\left|1 - \tau(p)^2\right| \left|\Q(h(p) + \delta \eta(p))\right|\\
 &\leq& \left|1 - \tau(p)^2\right|\frac{k _0}{2}\D( h(p)+\delta \eta(p))\\
  &\leq&\frac{3}{2} k_0 C  \left(\frac{1}{3t} \delta \psi(p) \epsilon^2 + \frac{2}{9t^2} 2^{7/3} \delta^2 \psi^2(p) \epsilon^4\right)  .
\end{eqnarray*}
Finally, from these estimates we get that for $p \in P$
$$\E(g(p))- \E(h(p)) \leq - 2 \delta \psi(p) \epsilon^{1/2} +  C_1\delta \psi(p) \epsilon^2, $$
where $C_1$ is a constant depending only on $k_0$, $t$ and $R$. Hence choosing at the beginning of the proof $\epsilon>0$ sufficiently small such that $ -2 \epsilon^{1/2}+C_1 \epsilon^2 < - \epsilon^{1/2}$
we get that
$$\E(g(p))- \E(h(p)) \leq -  \delta \psi(p) \epsilon^{1/2}.$$
If $p \notin P$ we have that $\psi(p)=0$ and $\E(g(p))=\E(h(p))$. If $\bar p \in \overline{B_R}$ is such that $\E(g(\bar p))=F(g)$, we have 
$$\E(h(\bar p))\geq \E(g(\bar p)) \geq c,$$
and hence $\bar p \in P$ and $\psi(\bar p)=1$. Thus, we get that
$$\E(g(\bar p))- \E(h(\bar p)) \leq -  \delta \epsilon^{1/2}$$
and in particular
$$F(g) + \epsilon^{1/2} \delta \leq \E(h(\bar p)) \leq F(h),$$
so that $g\neq h$. But by definition of $g$ we have 
$$d(g,h) \leq \delta$$
and hence
$$F(g) + \epsilon^{1/2} d(g,h) \leq F(h),$$
which gives a contradiction. The proof is complete.
\qed

\begin{Proposition}\label{propconstrminmaxPS2}
Let $K \in C^1(\R^3)$ satisfying $(K_1)$ and $(K_2)$, let $t\in\R^+$ and $R>0$ be fixed, and let $c$, $c_0$ be the numbers defined, respectively, in \eqref{cminimaxlevel}, \eqref{czero}. Assume that $c>c_0$. Then, for every sequence $(f_{n})\subset \Phi$ such that $\sup_{p \in \overline B_R} \E(f_n(p)) \to c$ there exists another sequence $(u_{n})\subset M_{t}$ such that $\E(u^n) \to c$ and with the additional property that
$$
\Delta u_{n}-K(u_{n})(u_{n})_{x}\wedge (u_{n})_{y}+\lambda (u_{n})_{x}\wedge (u_{n})_{y}\to 0\quad\text{in }\hat{H}^{-1}
$$
for some $\lambda\in\R$. 
\end{Proposition}
\Proof
Let $(f_{n})\subset \Phi$ be such that $\sup_{p \in \overline B_R} \E(f_n(p)) \to c$. Then, according to Proposition \ref{propconstrminmaxPS} we find sequences $(\epsilon_{n})\subset(0,1)$, with $\epsilon_{n}\to 0$, and $(u_{n})\subset M_{t}$ such that 
\begin{gather*}
c-\epsilon_n \leq\E(u^n) \leq \sup_{p \in \overline B_R} \E(f_n(p))\\
|\E^\prime(u^n)[\varphi]| \leq 2 \epsilon_n^{1/2}\|\varphi\|~~\forall \varphi \in T_{u^n} M_t \cap (\R^3 + C_0^\infty(\R^2,\R^3)) .
\end{gather*}
Then, since $(\R^3 + C_0^\infty(\R^2,\R^3)) \cap T_{u^n} M_t$ is dense in  $T_{u^n} M_t$ (see Lemma \ref{lemmadensity}) we conclude that 
\begin{equation}
\label{eq:Ekeland-inequality-2}
|\E'(u^{n})[\varphi]|\le2\epsilon^{1/2}_{n}\|\varphi\|~~\forall
\varphi\in T_{u^{n}}M_{t}~\!.
\end{equation}
Now let $v^{n}\in\hat{H}^{1}$ be the Riesz representative of $\V'(u^{n})$. Set 
$$
\lambda_{n}=\frac{\E'(u^{n})[v^{n}]}{\|v^{n}\|^{2}}
$$
(notice that $\lambda_{n}$ is well defined because $v^{n}\in L^{\infty}$, see Lemma \ref{L:volume}). For every $\varphi\in\hat{H}^{1}\cap L^{\infty}$ the projection of $\varphi$ on $T_{u^{n}}M_{t}$ is given by
$$
\tilde\varphi=\varphi-\frac{\langle v^{n},\varphi\rangle}{\|v^{n}\|^{2}}v^{n}
$$
and, by (\ref{eq:Ekeland-inequality-2}), 
$$
|\E'(u^{n})[\varphi]-\lambda_{n}\V'(u^{n})[\varphi]|=|\E'(u^{n})[\tilde\varphi]|\le2\varepsilon_{n}^{1/2}\|\tilde\varphi\|\le2\varepsilon_{n}^{1/2}\|\varphi\|,
$$
and then, by density, $\E'(u^{n})-\lambda_{n}\V'(u^{n})\to 0$ in $\hat{H}^{-1}$. Now we show that the sequence $(\lambda_{n})$ is bounded. First of all we observe that the sequence $(\D(u^n))$ is bounded, because $\E(u^n) \to c$ and by Remark \ref{R:K-volume} we know that $\E$ is coercive with constants depending only on $k_0$ (see also \eqref{eq2quaterpropconstrPS}). Thus, by (\ref{eq:Wente-inequality}), we estimate 
\begin{equation}
\label{eq:vn-bounded}
\|\nabla v^{n}\|_{2}+\|v^{n}\|_{\infty}\le C_1\|\nabla u^{n}\|_{2}^{2}\le C_2,
\end{equation} 
for some positive constants $C_1$, $C_2$.
Then
\begin{equation}
\label{eq:lambdan-bdd-1}
|\E'(u^{n})[v^{n}]|\le\left|\int_{\R^{2}}(\nabla u^{n}\cdot\nabla v^{n}+K(u^{n})v^{n}\cdot u^{n}_{x}\wedge u^{n}_{y})\right|\le\|\nabla u^{n}\|_{2}\|\nabla v^{n}\|_{2}+\|K\|_{\infty}\|v^{n}\|_{\infty}\|\nabla u^{n}\|_{2}^{2}\le C.
\end{equation}
Moreover, keeping into account that $\int_{\R^{2}}v^{n}\mu^{2}=0$ and being $\D(u^n)$ bounded, we have that
\begin{equation}
\label{eq:lambdan-bdd-2}
|3t|=|\V'(u^{n})[u^{n}]|=|\langle v^{n},u^{n}\rangle|=\left|\int_{\R^{2}}\nabla v^{n}\cdot\nabla u^{n}\right|\le\|\nabla v^{n}\|_{2}\|\nabla u^{n}\|_{2}\le C\|\nabla v^{n}\|_{2}=C\|v^{n}\|~\!.
\end{equation}
Then (\ref{eq:lambdan-bdd-1}) and (\ref{eq:lambdan-bdd-2}) imply that $(\lambda_{n})$ is bounded, because $t\ne 0$. Hence, for a subsequence $\lambda_{n}\to\lambda\in\R$ and since $(v^{n})$ is bounded in $\hat{H}^{1}$ (use (\ref{eq:vn-bounded})), we conclude that $\E(u^{n})-\lambda\V'(u^{n})\to 0$ in $\hat{H}^{-1}$.~
\hfill$\square$

\section{Proof of Theorem \ref{mainteo}}
In view of Remark \ref{R:corresp} we consider the functional $\F_K(u)= \A(u) + \Q(u)$ on $\hat{H^1}$. 
Let $t>0$ and denote by $\mathrm{Crit}_{\F_{K}}(t)$ the set of constrained critical points of $\F_{K}$ at volume $t$, which we define as
\begin{equation}\label{defconstrcritpF}
\begin{split}
\mathrm{Crit}_{\F_{K}}(t):=\Big\{u \in M_t~|~ &u_x\wedge u_y \neq 0 \ \hbox{a.e and}\\ 
&\exists\lambda \in \R \ \hbox{s.t.}\ \frac{d}{ds}\F_K(u+s\varphi)\Big|_{s=0} = \lambda \frac{d}{ds}\V(u+s\varphi)\Big|_{s=0}\  \forall \varphi \in C^{\infty}_0(\R^2,\R^3) \Big\}.
\end{split}
\end{equation} 

We point out that if $u$ is of class $C^2$ and free of branch points (i.e. $u$ parametrizes an immersed surface) then,  since $\varphi$ has compact support, we have
\begin{equation}\label{firstvararea}
\frac{d}{ds}\A(u+s\varphi )\Big|_{s=0} = -2 \int_{\R^2}H(u) \nu \cdot \varphi  |u_x\wedge u_y|,
\end{equation}
where $H$ is the mean curvature of $u$, $\nu=\frac{u_x\wedge u_y}{|u_x\wedge u_y|}$ is the Gauss map (see \cite{DHS}, Sect. 2.1, (7) and (8)). 

In general, if $u$ is smooth but not immersed then we can consider only variations $\varphi$ which have compact support in the set of regular points. Nevertheless, if $H$ is a prescribed function of class $C^{1,\alpha}$, then any $H$-bubble, nemaly any non constant (weak) solution $u\in\hat{H}^{1}$ of $\nabla u=2H(u)u_{x}\wedge u_{y}$ on $\R^{2}$, is in fact smooth, more precisely, of class $C^{3,\alpha}$, in view of well known results (see \cite{DHS}, Sect. 5.1, Theorem 1). Hence, the right-hand side of \eqref{firstvararea} can be continuously extended to variations $\varphi \in C^{\infty}_0(\R^2,\R^3)$. Therefore we can take  \eqref{firstvararea} as a definition of $\frac{d}{ds}\A(u+s\varphi )\Big|_{s=0}$ when $u$ is a $H$-bubble of class $C^{3,\alpha}$ (see also \cite{DHS}, Sect. 5.3).

Before proving Theorem \ref{mainteo} we need the following preliminary lemma.
 
\begin{Lemma}\label{techlemmacritpoint}
Let $K\in C^{1,\alpha}(\R^{3})$ satisfy $(K_{1})$ and $(K_{2})$. Then for any fixed $t>0$ it holds that $$\mathrm{Crit}_{\E}(t) \subset \mathrm{Crit}_{\F_K}(t).$$
\end{Lemma}
\Proof
If $u \in \mathrm{Crit}_{\E}(t)$, then by definition $u$ is a weak solution of 
$$
\Delta u=(K(u)-\lambda)u_{x}\wedge u_{y}\quad\text{on}\quad\R^{2},
$$
for some $\lambda \in \R$ and, by Lemma \ref{L:bubble-bdd}, $u$ is of class $C^{2,\alpha}$ as a map on $\S^2$ and satisfies the conformality relations $u_x \cdot u_y =0=|u_x|^2-|u_y|^2$ (see \cite{CaMu11}, Remark 2.5). Moreover, since we are assuming $K \in C^{1,\alpha}$, by well known regularity results (see Sect. 2.3, \cite{DHS2}), we get that $u$ is of class $C^{3,\alpha}$. Hence $u$ describes a closed parametric surface of mean curvature $\frac{1}{2}(K(u)-\lambda)$ in the set of regular points. 
Concerning the set of branch points of $u$ (i.e. points where $\nabla u =0$), we point out that it is at most finite (see \cite{GulOssRoy73} or \cite{DHS}, Sect. 5.1, \cite{DHS2}, Sect. 2.10), and in particular it holds that $u_x\wedge u_y \neq 0$ a.e. in $\R^2$. Since u is a $(K-\lambda)$-bubble of class $C^{3,\alpha}$, by \eqref{firstvararea} 
\begin{equation}\label{eqfvararea}
\frac{d}{ds}\A(u+s\varphi )\Big|_{s=0} = -2 \int_{\R^2} \frac{1}{2}(K(u)-\lambda) \nu \cdot \varphi |u_x\wedge u_y| = - \int_{\R^2} (K(u)-\lambda) \varphi \cdot u_x \wedge u_y,
\end{equation}
for any $\varphi \in C_0^\infty(\R^2,\R^3)$, where $\nu$ is the extension of the Gauss map (see \cite{DHS}, Sect. 5.1).
Now, from \eqref{eqfvararea} and Lemma \ref{L:K-volume} we get that for any $\varphi \in C_0^\infty(\R^2,\R^3)$
$$\frac{d}{ds}\F_K(u+s\varphi )\Big|_{s=0}=- \int_{\R^2} (K(u)-\lambda) \varphi \cdot u_x \wedge u_y + \int K(u) \varphi \cdot u_x \wedge u_y.$$
Moreover, by Lemma \ref{L:volume}, we have  $$\frac{d}{ds} \V(u+s\varphi )\Big|_{s=0} = \int_{\R^2} \varphi \cdot u_x \wedge u_y.$$
Hence, it immediately follows that
$$\frac{d}{ds}\F_K(u+s\varphi )\Big|_{s=0} = \lambda \frac{d}{ds} \V(u+s\varphi )\Big|_{s=0},$$
 for any $\varphi \in C_0^\infty(\R^2,\R^3)$, which means that $u \in \mathrm{Crit}_{\F_K}(t)$ (see \eqref{defconstrcritpF}).  
The proof is complete.~ 
\qed

Now we can prove Theorem  \ref{mainteo}.\\

\Proof
Assume by contradiction that the thesis is false. Then, by Lemma \ref{techlemmacritpoint}, there exists $t_{0} \in (0,\bar t)$ such that
 $$\mathrm{Crit}_{\E}(t) = \varnothing~~\forall t \in (0,t_{0}]~\!.$$ Hence the assumptions of Proposition \ref{propcrucial} are satisfied, and so there exists $R>0$ such that
\begin{equation}\label{eq1maintheo}
S t^{2/3} < c_0 < c <2^{1/3} S t^{2/3}~~\forall t \in (0,t_{0})~\!.
\end{equation}
By Proposition \ref{propconstrminmaxPS2}, there exists a constrained Palais-Smale sequence $(u^n) \subset M_t$ at level $c$.  Since $\D(u^n)$ is uniformly bounded (see the proof of Proposition \ref{propconstrminmaxPS2}), then, by Lemma \ref{L:PS-decomposition} we deduce that $I=\varnothing$ and $c=\sum_{j \in J} \D(U_j)$. 

Now we observe that, up to changing the index set $J$ we can assume that the coefficients $k_j \in \N^+$ in \eqref{eqalgebraicvolum} are all identically $1$. 
In fact for any given $j \in J$ if $k_j>1$ then we can split $\D(U_j)$ as the sum of the area of $k_j$ spheres having the same area $4\pi \lambda^2$ and the same volume $ \frac{4}{3} \pi \lambda^3$. Hence, up to replacing $j$ with $k_j$ new indexes $\tilde j_1,\ldots,\tilde j_{k_j}$ and repeating this operation for all $j \in J$ (we recall that $J$ is finite), then, we get a new finite index set $\tilde J$ such that all the algebraic multiplicities of the spheres $U_{\tilde j}$ are identically $1$.

Hence, denoting by $|\tilde J|$ the cardinality of $\tilde J$, we have $$c=\sum_{j \in J} \D(U_j)= \sum_{\tilde j \in \tilde J } \D(U_{\tilde j})=\sum_{\tilde j \in \tilde J} S t_{\tilde j}^{2/3}=S\left(\frac{t}{|\tilde J|}\right)^{2/3}|\tilde J|=S |\tilde J|^{1/3} t^{2/3},$$
but this contradicts \eqref{eq1maintheo}, because $|\tilde J|$ is a positive integer. The proof is complete.
\qed

As a consequence of Theorem \ref{mainteo}, and arguing as in the proof of Theorem 3.15 in \cite{Cald2015}, we get an existence result for the H-bubble problem.
\begin{Theorem}
\label{mainteo2}
Let $K\in C^{1,\alpha}(\R^{3})$ satisfy $(K_{1})$ with \eqref{condkzero}, $(K_{2})$, and assume that $K>0$ on $\R^3$. Then there exists a sequence $(\lambda_{n})\subset\R$ with $|\lambda_{n}|\to\infty$ such that for every $n$ there exists a $(K-\lambda_{n})$-bubble.
\end{Theorem}

\bigskip

\noindent
\textbf{Acknowledgements.} Work partially supported by the PRIN-2012-74FYK7 Grant ``Variational and perturbative aspects of nonlinear differential problems'', by the project ERC Advanced Grant 2013 n. 339958 ``Complex Patterns for Strongly Interacting Dynamical Systems - COMPAT'', and by the Gruppo Nazionale per l'Analisi Mate\-ma\-tica, la Probabilit\`a e le loro Applicazioni (GNAMPA) of the Istituto Nazionale di Alta Mate\-ma\-tica (INdAM).

\end{document}